\documentclass{amsart}

\usepackage{amssymb}
\usepackage{amsfonts}
\usepackage{amsmath}
\usepackage{euscript}
\usepackage{enumerate}
\usepackage{graphics}
\usepackage[active]{srcltx}

\newtheorem*{Theorem1'}{Theorem 1'}

\theoremstyle{definition}

\theoremstyle{remark}

\numberwithin{equation}{section}

\renewcommand \i{{\mathfrak i}}

\begin{document}

\title{Elementary criteria for irreducibility of $f(X^n)$}

\author{Natalio H. Guersenzvaig}
\address{Universidad CAECE (retired), Buenos Aires, Argentina.}
\email{nguersenz@fibertel.com.ar}

\begin{abstract}
Very simple sufficient conditions for the irreducibility of $f(X^n)$
 over an arbitrary unique factorization domain $Z$ are established via a generalization of a well known theorem of A. Capelli.
\end{abstract}
\maketitle

\vspace{-0.6cm}
\section{\bf Introduction}
We fix throughout this work a unique factorization domain $Z$ with field of fractions $Q$. The group of units
of $Z$ will be denoted by $U$. 

\smallskip
 Let $f(X)$ be any polynomial in $Z[X]$ of positive degree that is irreducible in $Z[X]$.  Using
the well known Eisenstein's Criterion we can easily show that, in
some cases, $f(X^n)$ will also be irreducible in $Z[X]$ for any
positive integer~$n$. However, this is not true in general. For
example, $f(X)=X^3-X^2-2X -1$ is irreducible in $Z[X]$ for $Z\in
\{\mathbb Z, \,\mathbb Z_2, \,\mathbb Z_5\}$, while
$f(X^2)=(X^3-X^2-1)(X^3+X^2+1)$.

In the main result of this work we will establish sufficient conditions
 for irreducibility of $f(X^n)$ in $Z[X]$ for any integer $n> 1$ that, besides $U$,
  only  depend on~$n$, the degree of $f(X)$ and the leading and constant
 coefficients of $f(X)$. (These conditions can be easily checked if $Z$ is an effective unique factorization domain.)
Elementary necessary and sufficient conditions will also be given.
(The adjective ``elementary'' refers to the fact that these conditions will be stated without
using proper algebraic extensions of~$Q$.)
The cases $f(X)= X-a$ and $f(X)=aX^2+bX+c$ are considered in [1, pp. 63-74].
Related results for polynomials over finite fields can be found in [7, pp. 93-95].

Henceforth we will use, for $S\subseteq Q$ and $t\in \mathbb N$, the following notations:
$$S^*=S\setminus \{0\}, \qquad S^t = \{s^t: s\in S\},\qquad tS=\{ts: s\in S\}.$$

Given $a, b\in Z$ and $m, n\in \mathbb N$ we denote C$(m, a, b, n)$ the following condition: 

\smallskip
For each prime $p$ dividing $n$ and any unit $u$ in~$U$ at least one of the two following statements is true:

\noindent (A) \,$ua\notin Z^p$;

\noindent (B) \,(i) \,$(-1)^mub\notin Z^p$ \,and\, (ii) \,$ub\notin Z^2$, if $4|n$.

\medskip
Our main result is the following theorem (it will be proved, together with an equivalent dual version,
in the last section of this paper). 

\medskip
\noindent{\sc THEOREM 1.1.} {\it Let $n$ be any integer greater than $1$, and let $f(X)$ be an arbitrary
polynomial in $Z[X]$ of positive degree~$m$, leading coefficient~$a$ and nonzero constant term~$b$, that is
irreducible in $Z[X]$. Assume that at least one of the conditions $\text{\rm C}(m, a, b, n)$, $\text{\rm C}(m, b, a, n)$ holds. Then
$$\text{\qquad $f(X^n)$ is irreducible in $Z[X]$.}$$}

\section{\bf Basic facts}
In this section we  review  some results that we will use later without specific reference.
First we remind the reader that the characteristic of $Z$, say $\chi(Z)$, is the only nonnegative integer
  that satisfies the following two conditions:
$$\text{(a) $\chi(Z)\cdot 1=0$; \,\,\, (b) if $k\in \mathbb Z$ and $k\cdot 1=0$, then $\chi(Z)|k$.}$$
\indent The following fact is needed to prove Corollary 4.6 below:

$\bullet$\, either $\chi(Z)=0$, or $\chi(Z)=p$ is a prime number in which case we have \linebreak
\indent\,\,\,\,\,\,$(x_1+ \cdots +x_n)^{p}=x_1^{p} + \cdots + x_n^{p}$ for all  $n\in \mathbb N$ and any $x_1, \dots ,
x_n\in Z$.

\smallskip
We also recall that a nonzero polynomial $f(X)\in Z[X]\setminus U$
 is called~{\it reducible in $Z[X]$} if there exist nonzero polynomials $g(X)$, $h(X)$ in $Z[X]\setminus U$ such that
 $f(X)=g(X)h(X)$. Otherwise $f(X)$ is called {\it irreducible in $Z[X]$}. The {\it content} of $f(X)$, say $c(f)$, is the greatest common divisor of their coefficients (modulo units of~$Z$), and $f(X)$ is called {\it primitive} if $c(f)=1$.
 Replacing $Z$ by $Q$ in this
 definition yields (since in this case $U=Q^*$) that $f(X)$ is irreducible
in $Q[X]$ if and only if $f(X)$ has positive degree and there are no polynomials $g(X)$, $h(X)$ in $Q[X]$ of
positive degree such that $f(X)=g(X)h(X)$.

\smallskip
 It is also well known the following result (see, for example, \cite{6}):

$\bullet$ \,if $f(X)\in Z[X]$ has positive degree, then, $f(X)$ is irreducible in $Z[X]$ if
\linebreak \indent\,\,\,\,\,\,and only if $f(X)$ is primitive (in $Z[X]$) and irreducible in $Q[X]$.

\smallskip
\noindent As a consequence, when $f(X)\in Z[X]$ has positive degree and it is irreducible in $Z[X]$ we can replace
$Z[X]$ by $Q[X]$ without risk in any of the expressions, ``$f(X^n)$ is reducible in $Z[X]$'', ``$f(X^n)$ is
irreducible in $Z[X]$''. To simplify, in these situations we will write ``$f(X^n)$ is reducible'' or
``$f(X^n)$ is irreducible'', respectively. For the same reason, except where the contrary is explicitly stated,
the terms ``primitive'' and ``prime'' will be understood  to apply to the sets $Z[X]$ and $\mathbb N$, respectively.

For any positive integer $n$ let $\Phi_n(X)$ denote the cyclotomic polynomial of order $n$ over $Q$, i.e., the polynomial
\begin{center}$\Phi_n(X)=\prod_{1\le k\le \phi(n)}(X-w^k)$.\end{center}
where $w$ denotes an arbitrary primitive $nth$ root of unity (in some splitting field of $\Phi_n(X)$ over $Q$), and $\phi$ denotes the Euler function (i.e., $\phi(n)$ is the number of integers in the set $\{1, \dots , n\}$ that are relatively prime to $n$).
The following facts (see  \cite{8}) are related to Theorem 4.3, Theorem 4.4 and Corollary 5.2:

$\bullet$ \,$\Phi_n(X)\in Z[X]$;

$\bullet$ \,$\Phi_n(X)$ is irreducible if $\chi(Z)= 0$;

$\bullet$ $X^{n}-1=\prod_{\substack{d\in \mathbb N\\ d|n}}\Phi_d(X)$, if either $\chi(Z)= 0$, or $\chi(Z)\not|n$;

$\bullet$ \,If $\chi(Z)\ne 0$, $\chi(Z)\not\!|n$ and $\#(Z)\!=\!q$, then $\Phi_n(X)$ factors into the product  \linebreak \indent\,\,\,\,\,\,of $\phi(n)/d$ distinct monic irreducible \mbox{polynomials} in $Z[X]$ of degree $d$,  \linebreak \indent\,\,\,\,\,\,where $d$ is the least positive integer $d$ such that $q^{d}\equiv 1\pmod{n}$.

\smallskip
We will use matrices and determinants as well. With $M_m(Q)$, $|A|$
 and $\Delta_A(X)$ we will respectively denote the ring of square matrices of order $m$ with coefficients in $Q$,
 the determinant of $A\in M_m(Q)$ and the characteristic polynomial of~$A$.
 In particular we will consider a well known type of matrices associated to polynomials.

\eject
 Let $f(X)$ be  an arbitrary polynomial in $Q[X]$ of positive
degree~$m$, say $f(X)\!=\!\sum_{j=0}^m a_jX^j$, and let $f^*(X)$ denote the monic polynomial associate to $f(X)$,
that is, $f^*(X)=\sum_{j=0}^mc_jX^j$, where $c_j\!=\!a_j/a_m$ for j\,=\,0, 1, \dots , $m$. The {\it companion
matrix} of $f^*(X)$, say $C_{f^*}$, is the matrix in $M_m(Q)$ defined by
\begin{equation*}
C_{f^*} = \begin{bmatrix}
0&1&\hdots&0&0\\
0&0&\hdots&0&0\\
\vdots&\vdots&\ddots&\vdots&\vdots\\
0&0&\hdots&0&1\\
-c_0&-c_1&\hdots&-c_{m-2}&-c_{m-1}
\end{bmatrix}.
\end{equation*}

We will freely use the following properties of $C_{f^*}$:

\smallskip
 $\bullet$ $f^*(X)$ is both the minimum polynomial of $C_{f^*}$ over~$Q$
 and the character-\linebreak \indent \,\,\,\,\,\,istic polynomial of~$C_{f^*}$ (so $f(X)=
a_m|XI_m-C_{f^*}|$);

$\bullet$ If $f(X)$ is
irreducible, then the ring $Q[C_{f^*}]=\{h(C_{f^*}): h(X)\in Q[X]\}$ is \linebreak \indent\,\,\,\,\,\,an extension field of $Q$ of degree $m$ with $f(C_{f^*})=a_mf^*(C_{f^*})=O$. \linebreak \indent\,\,\,\,\,\,(In this situation, as usual, we will write $Q(C_{f^*})$ instead of $Q[C_{f^*}]$.)

\section{\bf Preliminary results}

\smallskip
Our irreducibility criteria strongly depend on two
  beautiful theorems of A. Capelli (which are included in the author's Ph.~D. Thesis, Melbourne University, 1955).
 The first one gives non-elementary
necessary and sufficient conditions for irreducibility of $f(g(X))$ in $Q[X]$. For the sake of completeness we provide a simple proof of this result\footnote{\,Due to F. Szechtman, University of Regina, Saskatchewan, Canada}.

\medskip
\noindent{\sc CAPELLI'S THEOREM 1.} {\it Let $f(X)$, $g(X)$ be arbitrary polynomials of $Q[X]$ of positive degree. Let $F$ be any
splitting field of $f(X)$ over $Q$, and let $\alpha$ be any root of $f(X)$ in $F$. Then $f(g(X))$ is irreducible in $Q[X]$ if and
only if $f(X)$ is irreducible in $Q[X]$ and $g(X)-\alpha$ is irreducible in $Q(\alpha)[X]$.}
\begin{proof} We can assume without risk that $f(X)$ is irreducible in $Q[X]$. We also assume that $f(X)$ and
$g(X)$ have degrees $m$ and $n$, respectively, so $f(g(X))$ has degree $s=mn$. let $K$ be any splitting field of
$f(g(X))$ over~$F$. Letting $\alpha_1=\alpha$, we have in $K[X]$ the factorization
$$f(X) =a(X-\alpha_1)\dots (X-\alpha_m),$$
and hence, the factorization
$$f(g(X))= a(g(X)-\alpha_1)\dots (g(X)-\alpha_n)= b(X-\beta_1)\dots (X-\beta_s)$$
for certain nonzero $a, b\in Q$. Given any $i$ with $1\le i\le m$, there must exist at least one $j$, $1\le j\le s$, such
that $g(\beta_j)=\alpha_i$, for otherwise all $X-\beta_j$ would be relatively prime to $g(X)-\alpha_i$, and therefore
$f(g(X))$ would be relatively prime to $g(X)-\alpha_i$, a contradiction. This proves that $\alpha$ is of the form
$g(\beta)$, where $\beta$ is a root of $f(g(X))$ in~$K$. Then we have,
$$[Q(\beta): \,Q]=[Q(\beta): \,Q(\alpha)][Q(\alpha): \,Q]=[Q(\beta): \,Q(\alpha)]m.$$
\eject
\noindent Now, $f(g(X)$ is irreducible in $Q[X]$ if and only if $[Q(\beta): \,Q]=mn$, that is, if and only if $[Q(\beta):
\,Q(\alpha)]=n$. Therefore, since  $g(X)-\alpha$ has degree $n$ and annihilates $\beta$, $f(g(X))$ is irreducible in
$Q[X]$ if and only if $g(X)-\alpha$ is irreducible in $Q(\alpha)[X]$. This completes the proof.\end{proof}

\smallskip
 The second  establishes simple conditions for reducibility of $X^n-a$ in $Q[X]$ (see \cite{2} and \cite{7}, Theorem 9.1).

\medskip
 \noindent{\sc CAPELLI'S THEOREM 2.} {\it Let $a$ be any nonzero element of $Q$, and let~$n$ be any integer greater
 than 1.  Then $X^n -a$ is reducible in $Q[X]$ if and only if either} (i) {\it $a= c^t$ for some
  $c\in Q$ and $t|n$ with $t> 1$, or} (ii) {\it $4|n$ and $a= -4c^4$ for some $c\in Q$.}

\section{\bf Necessary and sufficient conditions}

\smallskip
In order to prove the main theorem of this section we first establish a result involving primitive polynomials
 that is interesting in its own right.

\medskip
\noindent{\sc LEMMA 4.1.} {\it Suppose that $P(X)=\sum_{k=0}^ma_kX^k$ is a primitive polynomial of $Z[X]$ of degree
$m$ with $a_0\ne 0$. Let $L(X)$ be any monic polynomial in $Z[X]$ of positive degree $n$ and nonzero roots
$\lambda_1, \dots , \lambda_n$ in some extension field of $Q$  (counting multiplicities). In addition suppose that
the constant coefficient of $L(X)$, say $c_0$, is relatively prime to $a_0$. Then
$$\text{$\prod_{j=1}^nP(\lambda_jX)$ is a primitive polynomial of $Z[X]$.}$$}

\vspace{-0.5cm}
\begin{proof} Since $L(X)$ is a monic polynomial of $Z[X]$, from the well known Fundamental Theorem on
Symmetric Polynomials it follows that $\prod_{j=1}^nP(\lambda_jX)$ is a polynomial in $Z[X]$, say
$$P^*(X)=\sum_{j=0}^{mn}a^*_jX^j.$$
\indent Looking for a contradiction we suppose $c(P^*)\ne 1$. Let $q$ be any prime of~$Z$ that divides $c(P^*)$.
As $c(P^*)$ divides $a^*_0=a^n_0$, $q|a_0$ and $q\not|c_0$. Thus, since $P(X)$ is primitive, there is a
positive integer $k$, $k\le m$, such that $q|a_j$ for $0\le j <k$ and $q\not|a_k$. Realizing the product
$\prod_{j=1}^{n}P(\lambda_jX)$ we get
$$a^*_{nk} = a_k^n\lambda^k_1\cdots \lambda^k_n \,+\,
\sum_{\substack{i_1+ \cdots + i_{n}=nk\\i_j\,\ge\, 0, \, j=1, \,\dots , \,n\\(i_1, \dots , \,i_{n})
 \ne(k, \dots , \,k)}}\!a_{i_1}\cdots a_{i_{n}}\lambda_1^{i_1}\cdots \lambda_n^{i_n}.$$

Notice that in each summand $a_{i_1}\cdots a_{i_{n}}\lambda_1^{i_1}\cdots \lambda_n^{i_n}$ we have $i_j<k$ for at
least one $j$, which makes each such summand a multiple of $q$. But $a^*_{nk}$ is also a multiple of~$q$. This
contradicts the fact that $a_k^n\lambda^k_1\cdots \lambda^k_n =(-1)^{nk}c_0^ka_k^{n}$ is not divisible by $q$.
\end{proof}

\medskip
 In addition we will use the following result,
which is an immediate consequence of two well known identities
(see, for example, (3.1.1)-(3.1.4) in \cite{3} and (22)-(25) in \cite{5}).

\eject
\noindent{\sc LEMMA 4.2.} {\it Let $F$ be an arbitrary field and let $p$ be any prime number.
Let $\Psi(Y)=Y^p-1\in F[Y]$ and let $w$ be any primitive $pth$ root of unity.
Let $g(Y)=\sum_{j=0}^{p-1}c_jY^j\in F(Y)$. Then
$$\prod_{j=0}^{p-1}g(w^j)=|g(C_{\Psi})|=\left|\begin{matrix}
c_0&c_1&\hdots&c_{p-2}&c_{p-1}\\
c_{p-1}&c_0&\hdots&c_{p-3}&c_{p-2}\\
\vdots&\vdots&\ddots&\vdots&\vdots\\
c_1&c_2&\hdots&c_{p-1}&c_0
\end{matrix}\right|.$$}
\indent Now the following extension  of Capelli's Theorem 2 can be proved.

\medskip
 \noindent{\sc THEOREM 4.3.} {\it Let $n$ be any integer, $n>1$, and let $f(X)$ be any irreducible
polynomial in $Z[X]$ of positive degree $m$ and leading coefficient $a$. The 
following two statements are equivalent.}

 (a) $f(X^n)$ {\it is reducible.}

 (b)  {\it There exist a prime $p$ that divides $n$, a unit $u$ in $U$ with $ua\in Z^p$, and
 \linebreak\indent\,\,\,\,\,\,\,\,\,\,\,polynomials $S_0(X)$, $S_1(X)$, \dots , $S_{p-1}(X)$ in $Z[X]$ such that either
  \indent\begin{align}\,\,\,
 (-1)^{m(p-1)}uf(X^{p})=\left|\begin{matrix}
S_{0}(X^p)\!\!&\!\!XS_1(X^p)& \!\hdots\!\!&\!X^{p-1}S_{p-1}(X^p)\\
X^{p-1}S_{p-1}(X^p)\!\!\!&\!\!S_0(X^p)& \!\hdots\! \!&\!X^{p-2}S_{p-2}(X^p)\\
\vdots \!\!\!&\!\!\vdots \!\!&\!\!\ddots \!\!&\!\vdots\\
XS_{1}(X^p)\!\!\!&\!\!X^{2}S_{2}(X^p)\!\!& \!\hdots \!\!&\!S_0(X^p)
 \end{matrix}\right|,
\end{align}
\indent\,\,\,\,\,\,\,\,\,or $4|n$ and
\begin{align}
uf(X^4)=\left|\begin{matrix}S_0(X^2)&XS_1(X^2)\\
XS_1(X^2)&S_0(X^2)\end{matrix}\right|.
 \end{align}}
\begin{proof} Assume (a). When $f(0)=0$, since $f(X)$ is irreducible, we have  $f(X)=aX$ with $a\in U$,
 so (1) follows with any prime $p$ that divides $n$, $u=a^{-1}$,
 $S_1(X)=1$  and $S_j(X)=0$ for $j=0, \dots , p-1$, $j\ne 1$. Therefore we may also
 assume $f(0)\ne 0$.

 Let $\alpha = C_{f^*}$. From Capelli's Theorem 1 it follows that $X^r-\alpha$
  is reducible in $Q(\alpha)[X]$. We first assume that condition (i) of Capelli's Theorem~1 holds. Therefore we have $\alpha=\gamma^t$  for some $\gamma\in Q(\alpha)$ and $t|n$, $t>1$.

Let $p$ be any prime that divides~$t$. Then we can write $\alpha = \beta^p$, where $\beta=\gamma^{t/p}\in
Q(\alpha)$. Hence, $X^p-\alpha$ is reducible in $Q(\alpha)[X]$, so $f(X^p)$ is reducible by Capelli's Theorem 1.

Let $\Psi(X)= X^p-1$ and let $w$ be an arbitrary primitive $pth$ root of unity.
 From $\Psi(X)= \prod_{j=0}^{p-1}(X-w^j)$ we get
 \begin{align*}
 X^p-\alpha &= X^p-\beta^p =\beta^p\Psi(\beta^{-1}X)= \prod_{j=0}^{p-1}(X-w^j\beta)\\
 &= w^{\frac{p(p-1)}2}\prod_{j=0}^{p-1}(w^{-j}X-\beta)=(-1)^{p-1}\prod_{j=0}^{p-1}(w^jX-\beta).
\end{align*}
Consequently, taking determinants on both sides, we obtain
 \begin{equation}
 f(X^{p}) =(-1)^{m(p-1)}a\prod_{j=0}^{p-1}\Delta_{\beta}(w^jX),
 \end{equation}
where $\Delta_{\beta}(X)=|XI_m-\beta|$, the characteristic polynomial of $\beta$, is a
 monic polynomial in $Q[X]$ of degree $m$. From unique factorization in $Z$ it follows that there exists $d\in Z$
  such that $P(X)=d\Delta_{\beta}(X)$ belongs to $Z[X]$ and  is
primitive. Since $P(X)$ has leading coefficient~$d$, letting $u=d^p/a$ (so $ua\in Z^p$) we can
rewrite (3) as follows:
 \begin{align}
 (-1)^{m(p-1)}uf(X^{p}) =\prod_{j=0}^{p-1}P(w^jX).
 \end{align}
The right hand side is a  primitive polynomial of $Z[X]$, by Lemma 4.1. Hence, since
$f(X^p)$ is also primitive (because $f(X)$ is), $u\in U$.

On the other hand, assuming $P(X)=\sum_{k=0}^ma_kX^k$ and expressing each index $k$ in the form $k=ip+j$ with $0\le
j< p$, we can write $a_kX^k= a_{ip+j}X^{ip}X^j$ for $k=0, \dots , m$. As a result, grouping the monomials associated
to each $X^j$ with $0\le j< p$, we obtain the polynomials
$$S_j(X)=\sum_{i\ge \,0}a_{ip+j}X^{i}\in Z[X],\,\,j=0, 1, \dots , p-1,$$
which satisfy
\begin{align}
P(X)=\sum_{j=0}^{p-1}X^jS_j(X^p).
\end{align}
Hence, since $C_{\Psi}^p$ is the identity matrix of order $p$, we get
 \begin{equation}
 P(XC_{\Psi})= \sum_{j= 0}^{p-1}X^jS_j(X^p)C_{\Psi}^j.
  \end{equation}
Thus (1) follows from the case $F= Q(X)$, $g(Y)=P(XY)$ of Lemma~4.2.

Now suppose that condition (ii) of Capelli's Theorem 2 applies. Thus $4\not|n$ and $\alpha =-4\gamma^4$ for some $\gamma\in Q(\alpha)$. It also follows that $f(X^2)$ is irreducible, because
otherwise case (i) applies with $p=2$. From the identity
$$X^4+ 4\gamma^4 = (X^2-2\gamma X + 2\gamma^2)(X^2+2\gamma X + 2\gamma^2)$$
it follows
$$f(X^4)=a|X^4I_m-\alpha| =a|X^2I_m-2\gamma X+ 2\gamma^2||X^2I_m + 2\gamma X+2\gamma^2|.$$
Hence, $f(X^4)$ is reducible in $Q[X]$. Consequently, condition~(i) of Capelli's Theorem 2 is satisfied with $n=2$ and $g(X)=f(X^2)$ instead of $f(X)$.
 By the first case for $g(X)$ there exist a unit $u$ in~$U$
  and $S_0(X)$, $S_1(X)$ in $Z[X]$ satisfying~(1) with $p=2$. Since this is the same as (2)
  (note that the degree of $g(X)$ is $2m$), we have completed the proof of (b).

Assume (b). Adding all other rows to the first row in the determinants of (1) and (2)
we get the row $[P(X) \hdots P(X)]$ (with $p=2$ in (2)), where $P(X)=\sum_{j=0}^{p-1}X^jS_j(X^p)$. Hence we can replace these determinants  by $P(X)P^*(X)$, where
\begin{align}P^*(X)=\left|\begin{matrix}
1\!\!&\!\!1&\!\!\hdots \!\!&\!1\\
X^{p-1}S_{p-1}(X^p)\!\!&\!\!S_0(X^p)&\hdots \!\!&\!X^{p-2}S_{p-2}(X^p)\\
\vdots \!\!&\!\!\vdots \!\!&\!\!\ddots \!\!&\!\vdots\\
XS_{1}(X^p)\!\!&\!\!X^{2}S_{2}(X^p)\!\!& \!\!\hdots \!\!&\!S_0(X^p)
 \end{matrix}\right|.\end{align}
Hence, for $i\in \{1, 2\}$, if (i) holds we have
\begin{equation}(-1)^{im(p-1)}uf(X^{ip}) = P(X)P^*(X),\end{equation}
i.e.,
 \begin{equation}uf(X^n)=
\begin{cases}
(-1)^mP(X^{n/p})P^*(X^{n/p})&\text{\!\!if (1) holds}\\
P(X^{n/4})P^*(X^{n/4})&\text{\!\!otherwise.}
\end{cases}\end{equation}
On the other hand, with $\Psi(X)$ and $w$ as defined above, from Lemma 4.2 we obtain,
\begin{equation}P(X)P^*(X)=|P(XC_{\Psi})|=\prod_{j=0}^{p-1}P(w^jX),\end{equation}
so that  $P(X)$ and $P^*(X)$ have positive degrees $im$ and $im(p-1)$, respectively.
This completes the proof of (a).\end{proof}

\medskip
Our next theorem provides complementary information about the polynomials $P(X)$ and $P^*(X)$
defined in Theorem 4.3.

\medskip
 \noindent{\sc THEOREM 4.4.} {\it Let $n$, $f(X)$, $\alpha$, $w$, $p$, $P(X)=\sum_{j=0}^{p-1}X^jS_j(X^p)$ and $P^*(X)$ be as defined in Theorem 4.3. Then}

(I) \,\,{\it $P(X)$ is irreducible.}

(II) {\it $P^*(X)$ is irreducible if and only if the cyclotomic polynomial
 $\Phi_p(X)=\indent \quad\,\,\,\,\sum_{j=0}^{p-1}X^j$
is irreducible and $Q(w)\cap Q(\alpha)=Q$.}
\begin{proof} We shall use (8). In order to prove (I) we suppose that $P(X)$ is reducible. Therefore, since $P(X)$ has degree $im$, $P(X)$ has a root,
say~$\gamma$, in some extension field of~$Q$ of degree $<im$ over~$Q$. Then $f(X^i)$ has also a root in
$Q(\gamma)$, namely~$\gamma^{p}$, contradicting the fact that $f(X^i)$ is irreducible.

Next we prove (II). This is clear if $p=2$, since $P(X)$ is irreducible and $P^*(X)=P(-X)$. Assume p is odd.

If $\Phi_p(X)$ is reducible, say
 $\Phi_p(X)=g(X)h(X)$, where both $g(X)$ and $h(X)$ are monic polynomials of positive degree  in $Z[X]$, from
$\Phi_p(X)=\prod_{j=1}^{p-1}(X-w^j)$ and Lemma 4.2 it follows the nontrivial factorization
$P^*(X)=|P(XC_{\Phi_p})|= |P(XC_{g})||P(XC_{h})|$.

\indent Assume $\Phi_p(X)$ is irreducible. Let $\delta=w\beta$, where $\beta$ is defined as in Theorem 4.3. Then, since $P(\beta)=0$,
$P^*(\delta)=\prod_{j=2}^{p}P(w^j\beta)=0$. Hence, $P^*(X)$ is
irreducible if and only if the minimum polynomial of $\delta$ over~$Q$ has degree $m(p-1)$.
 On the other hand we have $\delta\in Q(w, \beta)\subseteq Q(w, \alpha)$. Furthermore, since $\delta^p =
\beta^p=\alpha$, $\alpha\in Q(\delta)$ and  $Q(\alpha)=Q(\beta)$. Hence, $w=\beta^{-1}\delta\in Q(\delta)$. This
proves $Q(\delta) = Q(w, \alpha)$, and therefore that the minimum polynomial of~$\delta$ over~$Q$ has degree
$$[Q(\delta): Q]=[Q(w, \alpha): Q(w)][Q(w): Q]=[Q(w, \alpha): Q(w)](p-1).$$
From the well known theorem of natural irrationalities (see, for example, \cite{9}) we get $[Q(w, \alpha):
Q(w)]\!=\![Q(\alpha): Q(w)\cap Q(\alpha)]$, and hence
\begin{center}$[Q(\delta): Q]\!=\![Q(\alpha): Q(w)\cap Q(\alpha)](p-1)$.\end{center}
Now (II) follows immediately from $[Q(\alpha): Q]\!=\!m$. \end{proof}

\medskip
 \noindent{\it Remark.} In particular, $P^*(X)$ is  irreducible if $\Phi_p(X)$ is  irreducible and $\gcd(m, p-1)=1$.

\medskip
It should be noticed that in the course of the proof of Theorem 4.3 we have also proved, incidentally,
the following result.

\medskip
\noindent{\sc COROLLARY 4.5.} {\it Let $f(X)$ be any irreducible
polynomial in $Z[X]$ of positive degree. Let $n$ be any integer, $n>1$, and let $\sigma(n)$ be
the square-free part of~$n$. The three following statements are equivalent.}

\vspace{0.1cm}
\noindent (a) {\it $f(X^n)$ is reducible;}

\noindent (b) {\it either $f(X^{\sigma(n)})$ is reducible, or else $4|n$ and $f(X^4)$ is reducible;}

\noindent (c) {\it there exists a positive divisor of $r$, say
$t$, with $t$ prime or $t=4$, such that\linebreak \indent\quad\!\!\!$f(X^t)$ is reducible.

\smallskip
In particular, for any positive integer $s$ we have:}

\smallskip
\noindent (i) {\it $f(X^{2^s})$ is reducible if and only if either $f(X^2)$ is reducible, or else $s\ge 2$
and \linebreak \indent\quad\!\!\! $f(X^4)$ is reducible;}

\noindent (ii) {\it if $p$ is an odd prime, then
$$\text{$f(X^{p^s})$ is reducible if and only if $f(X^p)$ is reducible.}$$}
\indent For example, since $X^p$ can be replaced by $X$ in both sides of (1), we have:

\noindent (i) $f(X^{2^s})$ is reducible if and only if there exist $S_0(X)$, $S_1(X)$ in $Z[X]$ and
 \linebreak\indent\,\,\,$u\in U$ with $ua\in Z^2$ such that either
$$(-1)^muf(X)=S_0^2(X)-XS_1^2(X),$$
\indent\,\,\,\,or else
$$\text{$s\ge 2$ \,and \,$uf(X^2)=S_0^2(X)-XS_1^2(X)$.}$$

\noindent (ii) $f(X^{3^s})$ is reducible if and only if there exist $S_0(X)$, $S_1(X)$, $S_2(X)$ in
\linebreak\indent\,\,\,$Z[X]$ and $u\in U$ with $ua\in Z^3$ such that
$$uf(X) = S_0^3(X) +XS_1^3(X)+ X^2S_2^3(X)-3XS_0(X)S_1(X)S_2(X).$$
\indent On the other hand, from (II) of Theorem 4.4 it easily follows that (4) yields a factorization
of $f(X^{p})$ in $Z[X]$ into $p$ irreducible factors if and only if $w\in Z$.
From the case $w=1$, by using (i) and (ii) of Corollary 4.5 and the fact that
$$\text{$uf(X^2)\in Z^2[X]$ if and only if $uf(X)\in Z^2[X]$,}$$
the following result can also be easily derived.

\vspace{0.2cm}
\noindent{\sc COROLLARY 4.6.} {\it Assume $\chi(Z)=p$ is a prime number. Let $f(X)$ be any irreducible polynomial in $Z[X]$ of positive degree and let $s$ be any positive integer.}

 (a) {\it $f(X^{p})$ is reducible if and only if there exist
$u\in U$ and $P(X)\in Z[X]$, \indent\,\,\,\,\,\,\,\,\,$P(X)$ irreducible, such that
$$uf(X^p)=P^p(X);$$}
\indent  (b) {\it $f(X^{p^s})$ is reducible if and only if there exists $u\in U$ such that
$$uf(X) \in Z^p[X].$$}

\section{\bf Sufficient conditions}
First we derive Theorem 1.1 from Theorem 4.3.

\begin{proof} Assume $f(X)= \sum_{0\le k\le m}a_kX^k$ (so $a=a_m$, $b=a_0$). It is well known that either $f(X^n)=\sum_{0\le k\le m}a_kX^{nk}$ and  
$\tilde{f}(X^n)=X^{nm}f(1/X^{n})= \sum_{0\le k\le m}a_{m-k}X^{nk}$ 
are both irreducible, or both reducible in $Z[X]$. Thus we can assume that C$(m, a, b, n)$ holds. Looking for a contradiction suppose that $f(X^n)$ is reducible. In this situation, from Theorem~4.3 it follows that there exist a prime~$p$ that divides~$n$ ($p=2$ if (2) holds), a unit $u$ in $U$ with $ua\in Z^p$ (this contradicts (A)) and polynomials $P(X), P^*(X)$ in $Z[X]$ of positive degree that, in particular, satisfy $P(0)P^*(0) = (P(0))^p$  and (9). Therefore,  putting $X=0$ in both
sides of (9) we contradict (B) (i) if (1) holds, and (B)~(ii) otherwise. \end{proof}

As an immediate consequence of Theorem 1.1 we get the following result.

\medskip
\noindent{\sc COROLLARY 5.1.} {\it Let $f(X)$ be any non constant polynomial in $Z[X]$  with  leading coefficient~$a$ and nonzero constant term~$b$ that is
irreducible in $Z[X]$. If  $a\notin U$ or $b\notin U$, then the set of primes $p$ such that $f(X^p)$ is reducible is finite.}

\medskip
 When  $U$ is finite we may reformulate  Corollary 5.1  in the following way.

\medskip
\noindent{\sc COROLLARY 5.2.} {\it Let $f(X)$ be any non constant polynomial in $Z[X]$ with  leading coefficient~$a$ and nonzero constant term~$b$ that is
irreducible in $Z[X]$. Assume $U$ is finite. If $f(X)$ does not divides to any cyclotomic polynomial over $Q$, then the set of primes $p$ such that $f(X^p)$ is reducible is finite.}
\begin{proof} From Corollary 5.1 it will be sufficient to prove that both $a, b$ are in $U$ if and only if $f(X)$ is a divisor of some cyclotomic polynomial over $Q$.

First suppose $f(X)$ divides to some cyclotomic polynomial over $Q$, say $\Phi(X)$. Hence, since $\Phi(X)$ is a monic
polynomial in $Z[X]$ and its constant coefficient is a product of units in $U$, both $a$ and  $b$ belong to $U$.

Now suppose  $a\in U$ and $b\in U$. Since  the product of all roots of $f(X)$ (in some extension of $Q$) is equal to $\pm\,ba^{-1}\in U$), all roots of $f(X)$ in any extension of $Q$ are actually in $U$. Hence, for each root $\alpha$ of $f(X)$ we have $\alpha^s -1=0$, where $s= \#(U)$. Thus, since $f(X)$ is irreducible, $f(X)$ divides $X^s-1$. Therefore, since $X^s-1 =\prod_{\substack{d\in \mathbb N\\ d|s}} \Phi_d(X)$, $f(X)$ divides to some cyclotomic polynomial over $Q$. \end{proof}

\vspace{-0.1cm}
\noindent{\it Remark.} Hypothesis  $U$ is finite can not be suppressed in Corollary 5.2.  Because otherwise, considering for example,  $Z=Q= \mathbb R$ and the irreducible polynomial $f(X)=X^2+X+2$, we have that $f(X)$ does not divides to any cyclotomic polynomial over $Q$
(because the absolute values of the roots of $f(X)$ are distinct of $1$) and $f(X^p)$ is reducible for each prime $p$.

\vspace{0.2cm}
 At this point it should be noted that Theorem 1.1 essentially establishes that for a given positive
integer~$n$, if an arbitrary triple $(m, a, b)\in \mathbb N\times \mathbb Z^*\times \mathbb Z^*$
satisfies C$(m, a, b, n)$, then, for any $f(X)=aX^m +\cdots + b\in Z[X]$,
$$\text{$f(X^n)$ is irreducible if and only if $f(X)$ is irreducible.}$$

It is also of interest to determine, for a given irreducible polynomial $f(X)=aX^m + \cdots + b\in Z[X]$
of positive degree~$m$, an appropriate set of positive integers, say $\mathbb N(m, a, b)$, such that
 $f(X^n)$ is irreducible for each $r\in \mathbb N(m, a, b)$.

To illustrate the case that $f(X^n)$ is irreducible for all $r\in \mathbb N$ we consider Schur's polynomials, which
are defined for each positive integer $m$ by
$$f_m(X) = 1+ \frac {a_1}{1!} X+ \frac{a_2}{2!}X^2 + \cdots + \frac{a_{m-1}}{(m-1)!}X^{m-1} \pm\, \frac 1{m!}X^{m}\,\,
\,\,\,(\text{each $a_i\in \mathbb Z$)}.$$ It is well known that all these polynomials are irreducible in $\mathbb Q[X]$
(see [4, pp. 373-374]). It is clear that the polynomial 
$$m!f_m(X)= \pm\,X^m + ma_{m-1}X^{m-1} + \cdots + \frac{m!a_2}{2!}X^2 + \frac{m!a_1}{1!}X +m!\in \mathbb Z[X]$$
is primitive, so it is irreducible in  in $\mathbb Z[X]$. Assume $m\ge 2$. In some cases (for \linebreak example, when $m$
is prime) we get that $m!f_m(X^n)$ is irreducible for any $n\in \mathbb N$ from Eisenstein's Criterion, but in
general this does not happen (consider, for example, $m=2^n> 3$ and $a_{m-1}=1$). In any case we have $\pm
\,m!\not\in \mathbb Z^p$
 for each prime $p$, because the largest prime not exceeding~$m$ has such a property. Then, since
 condition (B) of C$(m, a, b, n)$ is always satisfied, m!$f_m(X^n)$ is irreducible
 (i.e., $f_m(X^n)$ is irreducible in $\mathbb Q[X]$) for any positive integer~$n$.

\smallskip
In order to include the precedent example in a more general result we assume that $a, \,b$ are arbitrary nonzero
elements of $Z$. First, we define the $(a, b)$-admissible primes. We shall say that a prime number $p$ is {\it $(a,
b)$-admissible} if there is no unit $u$ in $U$ such that both $ua$, $ub$ are in $Z^p$. Otherwise we shall say
that $p$ is {\it $(a, b)$-inadmissible}.

There is a simple procedure to determine the $(a, b)$-inadmissible primes. We first define the exponent of $(a,
b)$, say $e(a, b)$. Assume that $a$ has the factorization $a=u_ap_1^{\alpha_1}\cdots p_s^{\alpha_{s}}$ in $Z$,
where $u_a\in U$  and (in the case $a\not\in U$) $p_1, \dots , p_s$ are non-associate primes of $Z$ with positive
exponents $\alpha_1, \dots , \alpha_s$. Let $e(a)=0$ if $a=u_a$, and $e(a)=\gcd(\alpha_1, \dots , \alpha_s)$
otherwise. Assume a similar factorization for $b$, and let $e(a, b)=0$ if $e(a)=e(b)=0$ and $e(a, b)=\gcd(e(a),
e(b))$ otherwise. Then we can establish the following result.

\vspace{0.2cm}
 \noindent{\sc LEMMA 5.3.}  {\it Let $a$, $b$ be nonzero elements of $Z$ and let $p$ be a prime number. Then}
\begin{center}{\it $p$ is $(a, b)$-inadmissible if and only if $p|e(a, b)$ and $u_a\equiv u_b\!\!\pmod {U^p}$.}\end{center}
\begin{proof} To begin we express $a$ and $b$ in the form

\vspace{-0.2cm}
$$a =u_aa_0^{e(a)},\quad b=u_bb_0^{e(b)},$$ 

\vspace{-0.1cm}
\noindent where each one of $a_0$, $b_0$ is either equal $1$, or a product of non-associate prime-powers of $Z$.
Assume  $p|e(a, b)$ and $u_a^{-1}u_b\in U^p$, say $u_a^{-1}u_b=u_0^p$. Letting $e=e(a, b)$ we can write
$$a=u_a\alpha^e,\quad  b=u_b\beta^e,$$
where $\alpha=a_0^{e(a)/e}$, $\beta=b_0^{e(b)/e}$. Hence,
$$u_a^{e-1}a= (u_a\alpha)^e,\quad u_a^{e-1}b= (u_a^{-1}u_b)(u_a\beta)^e.$$
Thus $p$ is $(a, b)$-inadmissible, because
$$\text{$u_a^{e-1}a= ((u_a\alpha)^{e/p})^p$ \,and \,$u_a^{e-1}b=
(u_0(u_a\beta)^{e/p})^p$.}$$

Now assume that $p$ is $(a, b)$-inadmissible. Therefore, there exist $u\in U$ and $\alpha, \,\beta\in Z$ such that
$ua=\alpha^p$, $ub=\beta^p$. Proceeding as previously with $a$ and $b$, we can write
$\alpha=u_{\alpha}\alpha_0^{e(\alpha)}$, $\beta=u_{\beta}\beta_0^{e(\beta)}$, whence
$$uu_aa_0^{e(a)}=u_{\alpha}^p\alpha_0^{pe(\alpha)}, \quad uu_bb_0^{e(b)}=u_{\beta}^p\beta_0^{pe(\beta)}.$$
Hence, from the unique factorization property of $Z$, it follows $pe(\alpha)=e(a)$, $pe(\beta)=e(b)$ and both $uu_a,
\,uu_b\in U^p$. Thus,  $p|e(a, b)$ and $u_a^{-1}u_b\in U^p$.\end{proof}

 Next we define the $(a, b)$-admissible odd integers. For convenience we agree that $1$ is $(a,
b)$-admissible. Let $\mathbb N_{\text{o}}$ denote the set of odd positive integers. We shall say that $r\in
\mathbb N_{\text{o}}$ is $(a, b)$-admissible if each one of their prime divisors is $(a, b)$-admissible. Otherwise we
shall say that $r$ is $(a, b)$-inadmissible.

Let $\mathbb N_{\text{o}}(a, b)$ denote the set of $(a, b)$-admissible odd integers. The set $\mathbb N(m, a, b)$
of  {\it $(m, a, b)$-admissible} integers is defined then as follows:
$$\mathbb N(m, a, b)\! = \!\begin{cases}
\mathbb N_{\text{o}}(a, b) & \text{\,if $2$ is $(a, (-1)^mb)$-inadmissible,}\\
\mathbb N_{\text{o}}(a, b)\cup 2\mathbb N_{\text{o}}(a, b) & \text{\,if $2$ is both $(a, (-1)^mb)$-admissible}\\
&\text{\,and $(a, b)$-inadmissible,}\\
\cup_{k=0}^{\infty}2^k\mathbb N_{\text{o}}(a, b) &\text{\,if $2$ is both $(a, (-1)^mb)$-admissible}\\
&\text{\,and $(a, b)$-admissible}.\\
\end{cases}$$

Let $n$ be any integer greater than $1$. Writing $n=2^sq$, with $q$ odd and $s$ a nonnegative integer, we easily
get the following:

\smallskip
 \noindent (1) If $\mathbb N(m, a, b)=\mathbb N_{\text{o}}(a, b)$, then
$$r\in\mathbb N(m, a, b) \iff s=0 \text{ and } \text{C}(m, a, b, r);$$
\noindent (2) If $\mathbb N(m, a, b)=\mathbb N_{\text{o}}(a, b)\cup 2\mathbb N_{\text{o}}(a, b)$, then
$$r\in\mathbb N(m, a, b)\iff s\le 1 \text{ and } \text{C}(m, a, b, r);$$
\noindent (3) If $\mathbb N(m, a, b)=\cup_{k=0}^{\infty}2^k\mathbb N_{\text{o}}(a, b)$, then
$$r\in\mathbb N(m, a, b)\iff s\ge 0 \text{ and } \text{C}(m, a, b, r).$$
Hence,

\vspace{-0.5cm}
$$\text{$n\in \mathbb N(m, a, b)\iff \text{C}(m, a, b, n)$.}$$

\vspace{0.1cm}
Consequently we reformulate Theorem 1.1 as follows.

\medskip
 \noindent{\sc THEOREM 5.4.} {\it Let $n$ be any integer greater than $1$ and let $f(X)$
  be an irreducible polynomial in $Z[X]$ of positive
degree~$m$, leading coefficient~$a$ and nonzero constant term $b$.  Assume that at least one of the conditions $n\in \mathbb N(m, a, b)$,  $n\in \mathbb N(m, b, a)$ holds. Then
\begin{center}$f(X^n)$ is irreducible in $Z[X]$.\end{center}}

\vspace{0.1cm}
\indent Finally we use Lemma 5.3  to illustrate Theorem 5.4. Let $Z=\mathbb Z[i]$, where $i=\sqrt{-1}$,
 and assume $f(X)= X^m + \cdots + 8i$ is irreducible in $Z[X]$.  We have  $U=\{\pm 1, \pm i\}$, $a=1=u_a$, $e(a)=0$ and
$$\text{$b=-(2i)^3=-(1+i)^6$ with $u_b=-1$, $e(b)=6$.}$$
Then, since $e(a, b)=6$ and $u_a^{-1}u_b=-1\in U^p$ for each prime $p$, we have that~$2$ and $3$ are the unique
$(1, 8i)$-inadmissible primes. On the other hand, since  $(-1)^mu_a^{-1}u_b=(-1)^{m+1}\in U^2$, we also have
that $2$ is $(1, (-1)^m8i)$-inadmissible for all~$m$. Therefore,
\begin{center}$\mathbb N(m, 1, 8i)=
\mathbb N_{\text{o}}(1, 8i)=\{r\in \mathbb N_{\text{o}}: 3\not|r\}$,\end{center}
 which guarantees that $f(X^n)$ is irreducible for each positive integer $n$ that is relatively prime to $6$.



\end{document}